\documentclass[11pt,reqno]{amsart}

\usepackage{amssymb}
\usepackage{mathrsfs}
\usepackage{cite}
\usepackage{amsfonts}

\newcommand{\rmv}[1]{}

\def\wh{\widehat}

\def\pv#1{\ensuremath{{\bf#1}}}

\def\ilim{\varprojlim}

\def\inv{^{-1}}
\def\p{\varphi}

\def\pinv{{\p \inv}}

\def\R{\mathrel{{\mathscr R}}} 
\def\L{\mathrel{{\mathscr L}}} 

\def\e<{\leq _{E}}

\def\til#1{\ensuremath{\widetilde {#1}}}

\def\malce{\mathbin{\hbox{$\bigcirc$\rlap{\kern-8.3pt\raise0,50pt\hbox{$\mathtt{m}$}}}}}

\def\FP#1#2{\ensuremath{\widehat{F_{\mathbf #1}}(#2)}}

\def\1sk{^{(1)}}

\def\to{\rightarrow}

\def\stab#1{\mathrm{Stab}(#1)}
\def\hatexp#1#2{\widetilde{#2}}
\def\pl#1#2{\mathsf{PL}_{\pv {#1}}(#2)}

\def\Thmname{Theorem}
\def\Propname{Proposition}
\def\Lemmaname{Lemma}
\def\Definitionname{Definition}
\newtheorem{Thm}{\Thmname}[section]
\newtheorem{Prop}[Thm]{\Propname}
\newtheorem{Lemma}[Thm]{\Lemmaname}
\newtheorem{Def}[Thm]{\Definitionname}

\newtheorem{Cor}[Thm]{Corollary}

\numberwithin{equation}{section}

\title{A Profinite Approach to Stable Pairs}

\author{Karsten Henckell\and John Rhodes\and Benjamin Steinberg}
\address{Department of Mathematics/Computer Science \\ New College of Florida
5800 Bay Shore Road Sarasota, Florida 34243-2109\\ \and
Department of Mathematics\\
University of California at Berkeley \\
Berkeley \\ CA 94720\\
USA\\ \and School of Mathematics and Statistics\\
Carleton University \\
1125 Colonel By Drive\\
Ottawa, Ontario  K1S 5B6 \\
Canada}
\thanks{The third author was supported in part by NSERC} \email{bsteinbg@math.carleton.ca}
\date{December 7, 2006}

\dedicatory{Dedicated to the memory of Bret Tilson}

\subjclass{20M07}

\begin{document}
\begin{abstract}
We give a short proof, using profinite techniques, that idempotent
pointlikes, stable pairs and triples are decidable for the
pseudovariety of aperiodic monoids.  Stable pairs are also described
for the pseudovariety of all finite monoids.
\end{abstract}
\maketitle

\section{Introduction}
In this paper we introduce a new combinatorial technique for working
with elements of free pro-\pv V monoids where \pv V is a
pseudovariety of monoids closed under Malcev product on the left by
the pseudovariety of aperiodic monoids.  The approach uses the
Henckell-Sch\"utzenberger expansion, and essentially allows one to
transfer arguments from Combinatorics on Words to the profinite
context.  Let us describe some of the applications. Detailed
definitions are given below.

If $\pv V$ is a pseudovariety of monoids, then a finite monoid $M$
belongs to the Malcev product $\pv V\malce \pv A$ if and only if the
maximal $\pv A$-idempotent pointlikes of $M$ belong to $\pv
V$~\cite{PinWeil,Slice2,qtheor}. If \pv V is a local pseudovariety
of monoids~\cite{Tilson}, then $M$ belongs to the semidirect product
$\pv V\ast \pv A$ if and only if, for each maximal $\pv A$-stable
pair $(Y,N)$ of $M$, the quotient of $N$ by the kernel of the action
of $N$ on $Y$ belongs to \pv V~\cite{Henckellstable}. Henckell
proved that $\pv A$-idempotent pointlikes and $\pv A$-stable pairs
are computable~\cite{Henckellidem,Henckellstable}. We give a much
easier proof of his results using profinite techniques. Also we
characterize the stable pairs for the pseudovariety $\pv M$ of all
finite monoids, giving a partial answer to a question raised
in~\cite{Almeidahyp}.  We also prove that $\pv A$-triples
(introduced below) are decidable.

The paper is organized as follows.  First we introduce stable pairs
and pointlike sets, and prove a standard compactness result.  Next,
we recall the definition of the Henckell-Sch\"utzenberger expansion.
We then describe stabilizers in certain free pro-\pv V semigroups
and introduce a discontinuous homomorphism. This leads to a proof of
Henckell's theorem on idempotent pointlikes.  As a warm-up, we
handle \pv M-stable pairs before turning to \pv A-stable pairs. The
final section concerns aperiodic triples, which we believe will play
a role in the solution to deciding membership in the complexity one
pseudovariety.

\section{Stable pairs and pointlikes}
If $X$ is a set, we use $X^*$ for the free monoid, $X^+$ for the
free semigroup  and $\wh{X^*}$ for the free profinite monoid
generated by $X$~\cite{Almeida:book}.  If \pv V is a pseudovariety
of monoids we use $\FP V X$ to denote the free pro-\pv V monoid
generated by $X$. If $M$ is a monoid generated by a set $X$, then
the image of an element $w$ of $X^*$ or $\wh{X^*}$ (or $\FP V X$ if
applicable) in $M$ is denoted $[w]_M$.  As a shorthand, if
$\gamma\in \wh{X^*}$, then the image of $\gamma$ in $\FP V X$ is
denoted $[\gamma]_{\pv V}$.

If $M$ and $N$ are $X$-generated monoids, we define the
\emph{canonical relational morphism} $\p:M\to N$ by $n\in m\p$ if
and only if there exists $w\in X^*$ such that $[w]_M =m$ and
$[w]_N=n$; this is equivalent to there existing $\alpha\in \wh{X^*}$
with $[\alpha]_M = m$ and $[\alpha]_N=n$. If $M$ is an $X$-generated
monoid and $\pv V$ is a pseudovariety, then the canonical relational
morphism $\rho_{\pv V}:M\to \FP V X$ is the relational morphism
given by $\alpha\in m\rho_{\pv V}$ if and only if there exists
$\alpha'\in \wh{X^*}$ with $[\alpha']_M = m$ and $[\alpha']_{\pv V}
= \alpha$.  Alternatively, $\alpha\in m\rho_{\pv V}$ if and only if
there is a sequence $w_n\in X^*$ such that $w_n\to \alpha$ and
$[w_n]=m$ for all $n$.    We remark that $m\rho_{\pv V}$ is a closed
subset of $\FP V X$ (see~\cite{Slice,Slice2,RhodesStein,qtheor}).

\begin{Def}[$\pv V$-pointlikes]
If $M$ is a finite monoid and $Z\subseteq M$ is a subset, then
$Z$ is said to be $\pv V$-pointlike, if for all relational
morphisms $\p:M\to N$ with $N\in \pv V$, there exists $n\in N$ such
that $Z\subseteq n\pinv$.
\end{Def}

The collection $\pl V M$ of $\pv V$-pointlikes of $M$ is a submonoid
of the power set $P(M)$, which is downwards closed in the order
$\subseteq$. The following fact about pointlike sets is well known.
Proofs can be found in~\cite{Henckell,Slice,Slice2}, for instance.

\begin{Prop}\label{liftandpushpoint}
Let $\pv V$ be a pseudovariety of monoids.   The map $M\mapsto \pl V
M$ is a functor preserving onto maps.  More precisely, if $\p:M\to
N$ is a homomorphism and $Z\in \pl V M$, then $Z\p\in \pl V N$.  If,
in addition, $\p$ is onto, then given $Z\in \pl V N$, there exists
$Z'\in \pl V M$ with $Z'\p = Z$.
\end{Prop}

So given a homomorphism $\p:M\to N$, there is an induced
homomorphism $\p_*:\pl V M\to \pl V N$ given by $Z\p_* = Z\p$ (the
direct image).




An element $Z\in \pl V M$ is called \pv V-\emph{idempotent
pointlike} if, for all relational morphisms $\p:M\to N$ with $N\in
\pv V$, there exists an idempotent $e\in N$ with $Z\subseteq
e\pinv$.  Notice that if $Z\in \pl V M$ and $Z=Z^2$, then $Z$ is
trivially idempotent pointlike since if $Z\subseteq n\pinv$, then
$Z\subseteq n^{\omega}\pinv$.  Also the set of \pv V-idempotent
pointlikes of $M$ form a downwards closed subset of $P(M)$.

Next we consider the notion of a \pv V-stable pair.  If $M$ is a
monoid and $s\in M$, then the stabilizer of $m$ is the
submonoid \[\stab m=\{m'\in M\mid mm'=m\}.\]

\begin{Def}[$\pv V$-stable pairs]
 Let
$M$ be a monoid. A pair $(Y,U)$ with $Y\subseteq M$ and $U\leq M$
(a submonoid) is called a $\pv V$-stable pair if, for all
relational morphisms $\p:M\to N$ with $N\in \pv V$, there exists
$n\in N$  such that $Y\subseteq n\pinv$ and $N\leq \stab n\pinv$.
\end{Def}

 If
we use the product ordering on pairs $(Y,N)$ with $Y$ a subset and
$N$ a submonoid of $M$, then the set of $\pv V$-stable pairs is
downwards closed.  Notice that to decide which pairs are \pv
V-stable, we just need to be able to compute all the \emph{maximal}
ones. Therefore, we focus our attention on these.  Observe that if
$(Y,U)$ is a stable pair, then so is $(YU,U)$.  Thus the maximal
stable pairs are transformation monoids.  It is straightforward to
verify that if $\pv V$ is a local pseudovariety of monoids in the
sense of Tilson~\cite{Tilson}, then $M\in \pv V\ast \pv W$ if and
only if, for each maximal $\pv W$-stable pair $(Y,U)$ of $M$, the
quotient of $U$ by the kernel of the action on $Y$ belongs to $\pv
V$, c.f.~\cite{Henckellstable,HMPR}.

Let us consider a more general notion.  A directed graph $\Gamma$
consists of a vertex set $V(\Gamma)$, an edge set $E(\Gamma)$ and
functions $\iota,\tau:E(\Gamma)\to V(\Gamma)$ selecting the initial
and terminal vertices $e$ of an edge, respectively.  We consider
only finite graphs. A \emph{labelling} of a graph $\Gamma$ over a
monoid $M$ is a function $\ell:V(\Gamma)\cup E(\Gamma)\to P(M)$. If
the image of $\ell$ is contained in $M$, we call $\ell$ a
\emph{singleton labelling}.  A singleton labelling $\ell$ is said to
\emph{commute} if, for each edge $e$, $e\iota\ell e\ell =
e\tau\ell$. If $\p:M\to N$ is a relational morphism, $\ell$ is a
labelling of $\Gamma$ over $M$ and $\ell'$ is a singleton labelling
of $\Gamma$ over $N$, then $\ell$ is said to be \emph{$\p$-related}
to $\ell'$ if $x\ell\subseteq x\ell'\pinv$ for all vertices and
edges $x$ of $\Gamma$. The following notion generalizes a notion of
Almeida~\cite{Almeidahyp}, which in turn generalizes a notion of
Ash~\cite{Ash}.

\begin{Def}[\pv V-inevitable graph]
Let $M$ be a finite monoid and \pv V a pseudovariety. A labelling
$\ell$ of a graph $\Gamma$ over $M$ is $\pv V$-inevitable if, for
all relational morphisms $\p:M\to N$ with $N\in \pv V$, there is a
singleton labelling $\ell'$ of $\Gamma$ over $N$ which commutes and
which is $\p$-related to $\ell$.
\end{Def}

For instance, $Z\subseteq M$ is $\pv V$-pointlike if and only if the
graph with a single vertex labelled by $Z$ is \pv V-inevitable. Let
$Y\subseteq M$ and $N\leq M$.  Let $\Gamma$ be a graph with one
vertex and $|N|$ loops. Then $(Y,N)$ is a \pv V-stable pair if and
only if the labelling of $\Gamma$ that assigns $Y$ to the vertex and
labels the edges by the elements of $N$ is \pv V-inevitable.
Conversely, a labelling of a graph with one vertex with  label $Y$
and that assigns singletons to the loops at the vertex is $\pv
V$-inevitable if and only if $(Y,\langle Z\rangle)$ is a \pv
V-stable pair where $Z$ is the set of labels of the edges. A
singleton labelling of a graph $\Gamma$ by $M$ is \pv V-inevitable
if and only if it is \pv V-inevitable in the sense of
Almeida~\cite{Almeidahyp}. Conversely, one can go from inevitable
labellings in our sense to that of Almeida by changing the graph.
For instance, $Z\subseteq M$ is \pv V-pointlike if and only if the
singleton labelling of a graph with two vertices and $|Z|$ directed
edges, where the initial vertex is labelled $1$, the edges are
labelled by the elements of $Z$ and the terminal vertex is labelled
by some element of $Z$, is \pv V-inevitable~\cite{Almeidahyp}.
Similarly, $(Y,N)$ is a \pv V-stable pair, if and only if the
singleton labelling of the graph $\Gamma$ with two vertices
$v_1,v_2$, $|Y|$ edges from $v_1$ to $v_2$ and $|N|$ loops from
$v_2$ to $v_2$ where $v_1$ is labelled by $1$, $v_2$ is labelled by
some element of $Y$, the $|Y|$ edges are labelled by the elements of
$Y$ and the $|N|$ loops are labelled by the elements of $N$ is \pv
V-inevitable. We leave the general construction to the reader.

Our notion has the advantage that it is closed downwards in the
partial order. That is, a labelling $\ell:V(\Gamma)\cup E(\Gamma)\to
P(M)$ can be viewed as an element of $P(M)^{V(\Gamma)\cup
E(\Gamma)}$.  If we order this set by the product ordering, then the
$\pv V$-inevitable elements form a down-set.  Decidability is then
reduced to calculating the \emph{maximal} elements.

The next two results give the relationship between the notions we
have been discussing and profinite techniques.  We include them for
completeness, and readers already conversant with this subject
should feel free to skip them.

\begin{Lemma}\label{canonical}
Let $M$ be a finite $X$-generated monoid and \pv V a pseudovariety
of monoids.  Let $\rho_{\pv V}:M\to \FP V X$ be the canonical
relational morphism.  Write $\FP V X= \ilim M_{\alpha}$ where the
$M_{\alpha}$ are $X$-generated monoids in \pv V.   Let
$\rho_{\alpha}:M\to M_{\alpha}$ be the canonical relational morphism
and \mbox{$\pi_{\alpha}:\FP V X\to M_{\alpha}$} the canonical
projection for each $\alpha$. Then:
\begin{enumerate}
\item If $C = \ilim C_{\alpha}\subseteq \FP V X$ is an inverse limit
  of subsets $C_{\alpha}\subseteq M_{\alpha}$ (with the induced
  inverse system), then $C\rho_{\pv V}\inv = \bigcap
  C_{\alpha}\rho_{\alpha}\inv$;
\item $\rho_{\pv V}\inv = \bigcap \pi_{\alpha}\rho_{\alpha}\inv$;
\item If $\gamma \in \FP V X$, then
 $\stab{\gamma}=\ilim \stab {\gamma\pi_{\alpha}}$.
\end{enumerate}
\end{Lemma}
\begin{proof}
Since $\rho_{\alpha} = \rho_{\pv V}\pi_{\alpha}$ and
$C\pi_{\alpha}\subseteq C_{\alpha}$, we have, for all $\alpha$,
\[C\rho_{\pv V}\inv \subseteq C\pi_{\alpha}\pi_{\alpha}\inv \rho_{\pv
V}\inv = C\pi_{\alpha}\rho_{\alpha}\inv \subseteq
C_{\alpha}\rho_{\alpha}\inv.\] For the converse, suppose $m\in
\bigcap
  C_{\alpha}\rho_{\alpha}\inv$.  Let $Y_{\alpha} = \{y\in
  C_{\alpha}\mid m\in y\rho_{\alpha}\inv\}$.  Then the $Y_{\alpha}$ are
  easily verified to form an inverse system.  By assumption on $m$,
  the $Y_{\alpha}$ are non-empty finite sets.  Hence $\emptyset \neq
  \ilim Y_{\alpha}\subseteq \ilim C_{\alpha}=C$.   Now $(\ilim
  Y_{\alpha})\pi_{\beta}\subseteq Y_{\beta}\subseteq m\rho_{\beta} = m\rho_{\pv
    V}\pi_{\beta}$, for all $\beta$, and hence $\ilim
  Y_{\alpha}\subseteq m\rho_{\pv V}$, since $m\rho_{\pv V}$ is
  closed~\cite{Slice,Slice2,RhodesStein,qtheor}.  This shows that $m\in
  C\rho_{\pv V}\inv$ and completes the proof of (1).

One deduces (2) from (1) by observing that if $\gamma\in \FP V X$,
then $\{\gamma\} = \ilim \{\gamma\pi_{\alpha}\}$.
Item (3) is clear from the description of $\ilim M_{\alpha}$ as a
subsemigroup of $\prod M_{\alpha}$ (see also~\cite[Proposition
9.6]{RhodesStein}).
\end{proof}

The following compactness result encompasses several well-known such
results~\cite{Almeidahyp,qtheor,Slice,Slice2,PinWeil}.

\begin{Thm}\label{profinitewitness}
Let $M$ be a finite $X$-generated monoid and \pv V a pseudovariety of
monoids.  Let $\rho_{\pv V}:M\to \FP V X$ be the canonical relational
morphism.  Then:
\begin{enumerate}
\item A subset $Z\subseteq M$ is \pv V-pointlike if and only if there
  exists $\alpha\in \FP V X$ with $Z\subseteq \alpha\rho_{\pv V}\inv$;
\item A subset $Z\subseteq M$ is \pv V-idempotent pointlike if and
  only if there
  exists an idempotent $\alpha\in \FP V X$ with $Z\subseteq \alpha\rho_{\pv
    V}\inv$;
\item $(Y,N)$ is a $\pv V$-stable pair for $M$ if and only if there
  exists $\alpha\in \FP V X$ with $Y\subseteq \alpha\rho_{\pv V}\inv$
  and $N\leq \stab {\alpha}\rho_{\pv V}\inv$;
\item A labelling $\ell$ of a graph $\Gamma$ is $\pv V$-inevitable if
  and only if there is a singleton labelling of $\Gamma$ over $\FP V X$ that is
  $\rho_{\pv V}$-related to $\ell$ and which commutes.
\end{enumerate}
\end{Thm}
\begin{proof}
We prove (3) and (4).  A proof of (1) and (2) can be found
in~\cite{Slice2} (alternatively, (2) follows from (4)).

For (3),  let $Y\subseteq M$ and $N\leq M$. If $\p:M\to T$ is a
relational morphism, there is always an $X$-generated submonoid $T'$
of $T$ and a canonical relational morphism $\psi:M\to T'$ such that
$\psi\subseteq \p$ (as relations)~\cite{qtheor}.  So we may take all
the relational morphisms in the definition of a $\pv V$-stable pair
to be canonical relational morphisms of $X$-generated monoids.
Suppose that $\FP V X= \ilim M_{\alpha}$ where the $M_{\alpha}$ run
over all $X$-generated monoids in $\pv V$. Let
\mbox{$\rho_{\alpha}:M\to M_{\alpha}$} and \mbox{$\rho_{\pv V}:M\to
\FP V X$} be the canonical relational morphisms and denote by
\mbox{$\pi_{\alpha}:\FP V X\to M_{\alpha}$} the canonical
projection.  Set
\begin{align*}
C_{\alpha} &= \{m\in M_{\alpha}\mid Y\subseteq m\rho_{\alpha}\inv,\
N\leq \stab m\rho_{\alpha}\inv\} \\ C &= \{\gamma\in \FP V X\mid
Y\subseteq \gamma\rho_{\pv V}\inv,\ N\leq \stab{\gamma}\rho_{\pv
V}\inv\}.
\end{align*}
Then the $C_{\alpha}$ are easily verified to form an inverse system.
  We claim that $C =\ilim C_{\alpha}$.  Since an
inverse limit of finite sets is non-empty if and only if each of the
finite sets is non-empty, this will yield (3). Indeed, applying
Lemma~\ref{canonical} we see that, for $\gamma\in \FP V X$, the
equalities
\begin{align*}
\gamma\rho_{\pv V}\inv &= \bigcap
\gamma\pi_{\alpha}\rho_{\alpha}\inv\\
 \stab{\gamma}\rho_{\pv V}\inv & = \bigcap
  \stab{\gamma\pi_{\alpha}}\rho_{\alpha}\inv
\end{align*}
hold. Thus $\gamma\in C$ if and only if $Y\subseteq \gamma\rho_{\pv
V}\inv$, $N\leq \stab{\gamma}\rho_{\pv V}\inv$, if and only if
$Y\subseteq \gamma\pi_{\alpha}\rho_{\alpha}\inv$, $N\leq
\stab{\gamma\pi_{\alpha}}\rho_{\alpha}\inv$ all $\alpha$, if and
only if $\gamma\in \ilim C_{\alpha}$, as required.


For (4), let $\Gamma$ be a graph. If $N$ is a monoid, we use
$N^{\Gamma}$ as a shorthand for $N^{V(\Gamma)\cup E(\Gamma)}$.  As
before, we need only consider canonical relational morphisms of
$X$-generated monoids when considering $\pv V$-inevitability.

Consider a labelling $\ell\in P(M)^{\Gamma}$.  Write $\FP V X= \ilim
M_{\alpha}$ where the $M_{\alpha}$ run over all $X$-generated
monoids in \pv V. Let $\rho_{\alpha}:M\to M_{\alpha}$ and
\mbox{$\rho_{\pv V}:M\to \FP V X$} be the canonical relational
morphisms. Let $C_{\alpha}(\Gamma)\subseteq M_{\alpha}^{\Gamma}$ be
the set of all commuting singleton labellings of $\Gamma$ that are
$\rho_{\alpha}$-related to $\ell$. Then the $C_{\alpha}(\Gamma)$
form an inverse system.   Indeed, if
$\pi_{\alpha,\beta}:M_{\alpha}\to M_{\beta}$ is the canonical
projection, then the image under $\pi_{\alpha,\beta}$ of a commuting
singleton labelling of $M_{\alpha}$ clearly commutes and also
$\rho_{\alpha}\inv \subseteq \pi_{\alpha,\beta}\rho_{\beta}\inv$ so
$\rho_{\alpha}$-related labellings to $\ell$ are sent to
$\rho_{\beta}$-related labellings.

Let $C(\Gamma)\subseteq \FP V X^{\Gamma}$ be the set of all
commuting singleton labellings of $\Gamma$ that are $\rho_{\pv
V}$-related to $\ell$.  Then $C(\Gamma)$ is a closed subset of the
profinite monoid $\FP V X^{\Gamma}$ and, in fact, $C(\Gamma) = \ilim
C_{\alpha}(\Gamma)$.  Indeed, writing $\pi_{\alpha}:\FP V X\to
M_{\alpha}$ for the canonical projection, we have that $\rho_{\pv
V}\inv = \bigcap \pi_{\alpha}\rho_{\alpha}\inv$ by
Lemma~\ref{canonical} and a labelling $\ell'\in \FP V X^{\Gamma}$
commutes if and only if all its images in the $M_{\alpha}$ commute
(viewing $\FP V X$ as a submonoid of $\prod M_{\alpha})$.  Since the
inverse limit of an inverse system of finite sets is non-empty if
and only if each of the sets is non-empty, we conclude that $\ell$
is $\pv V$-inevitable if and only if $C(\Gamma)\neq\emptyset$.  This
completes the proof.
\end{proof}

\section{The Henckell-Sch\"utzeneberger Expansion}
Our key tool for understanding stable pairs and related notions is the
Henckell-Sch\"utzenberger expansion.  Further applications of this
expansion can be found in~\cite{projective}.
Recall that if $M$ and $N$ are monoids, then their
\emph{Sch\"utzenberger product}~\cite{Eilenberg,BR--exp} is the
monoid
\[M\lozenge N = \begin{bmatrix} M & P(M\times N) \\ 0 &
N\end{bmatrix}\] with multiplication given by
\[\begin{bmatrix} m & U
 \\ 0 & n\end{bmatrix}\begin{bmatrix} m' & U'
 \\ 0 & n'\end{bmatrix}=\begin{bmatrix} mm' & mU'+Un'
 \\ 0 & nn'\end{bmatrix}\]
where addition is union and where $P(M\times N)$ is viewed as an
$M$-$N$-bimodule in the obvious way.

If $M$ is an $X$-generated monoid, then the \emph{Henckell-Sch\"utzenberger
  expansion} $\hatexp 2 M$ is the submonoid of $M\lozenge M$ generated
by matrices of the form \[\begin{bmatrix} x & (1,x)+(x,1)\\ 0 &
  x\end{bmatrix}\] with $x\in X$.  So $\hatexp 2 M$ is an $X$-generated
monoid mapping naturally onto $M$ via the projection $\eta:\hatexp 2
M\to M$ to the diagonal.  Thus $\hatexp 2 M$ is an expansion
cut-to-generators in the sense of~\cite{BR--exp}.  Since $M\lozenge
M$ is really a double semidirect product of $M$ with $P(M\times M)$,
it follows $\eta$ is an $\pv {LSl}$-morphism~\cite{qtheor}, meaning
that the inverse image of each idempotent is locally a semilattice.
In particular $\eta$ is an aperiodic morphism (see
also~\cite{BR--exp,Eilenberg}).

Let $w\in X^+$.  By a \emph{cut} of $w$ we mean a pair $(u,v)\in
X^*\times X^*$ such that $w=uv$.  The set of cuts of $w$ will be
denoted $\vec c(w)$; we set $\vec c(\varepsilon)=\emptyset$. The
next proposition is well known~\cite{Eilenberg,BR--exp} and can be
proved by a simple induction on length.

\begin{Prop}\label{congruencehat}
Let $w\in X^*$ and $M$ an $X$-generated monoid.  Then \[[w]_{\hatexp 2
  M} = \begin{bmatrix} [w]_M & \sum_{(u,v)\in \vec c(w)}
  ([u]_M,[v]_M)\\ 0 & [w]_M\end{bmatrix}.\]  In particular, for
  $w,w'\in X^+$, the equality
$[w]_{\hatexp 2 M} = [w']_{\hatexp 2 M}$ holds if and only if, for
each factorization $w=uv$, there is a factorization $w' = u'v'$ such
that $[u]_M = [u']_M$ and $[v]_M = [v']_M$, and vice versa.
\end{Prop}

Henckell~\cite{Henckellstable} observed that stabilizers in $\hatexp 2
M$ enjoy a certain nice property.

\begin{Lemma}[Henckell]\label{criticallemmaforpairs}
Let $M$ be an $X$-generated monoid and let $w\in X^*$.
Then $\mathrm{Stab}([w]_{\hatexp 2 M})\eta$ is an $\L$-chain in the
monoid $\mathrm{Stab}([w]_M)$.
\end{Lemma}
\begin{proof}
If $w=\varepsilon$, there is nothing to prove, so assume $w\in X^+$.
Suppose that $u,v\in X^*$ with $[wu]_{\hatexp 2 M} = [w]_{\hatexp 2
M}=[wv]_{\hatexp 2 M}$. Then $wu=w_1w_2$ where $[w_1]_M = [w]_M$ and
$[w_2]_M=[v]_M$. There are two cases. Suppose first that $|w|\leq
|w_1|$.  Then $w_1=wx$ and $wxw_2=wu$, so $u=xw_2$.  Therefore,
$[w]_M=[w_1]_M = [w]_M[x]_M$, establishing that $[x]_M\in \stab
{[w]_M}$. In addition, $[u]_M = [x]_M[w_2]_M = [x]_M[v]_M$ and so
$[u]_M\leq_{\L} [v]_M$ in $\stab {[w]_M}$. If $|w|>|w_1|$, then $w =
w_1y$ and $w_1w_2 = w_1yu$, so $w_2 = yu$. A similar argument to the
above one then shows that $[y]_M\in \stab {[w]_M}$ and
$[v]_M\leq_{\L} [u]_M$ in $\stab {[w]_M}$. This completes the proof.
\end{proof}

\section{The structure of stabilizers and idempotent pointlikes}
We begin with some applications of the Henckell-Sch\"utzenberger
expansion to stabilizers and idempotent pointlikes.

\subsection{The structure of stabilizers}
Our first goal is to characterize stabilizers for free pro-$\pv V$
semigroups when $\pv V=\pv A\malce \pv V$, that is, $\pv V$ is
closed under the Henckell-Sch\"utzenberger expansion. The approach
is similar to the one taken in~\cite{RhodesStein} for related
results. A monoid $M$ will be called an \emph{internal $\L$-chain}
if the
 $\L$-classes of $M$ form a chain for the $\L$-ordering.  The reason
 the word internal is used is because if $M\leq N$, then $M$ can be an
 $\L$-chain in $N$ without being an internal $\L$-chain.

\begin{Thm}\label{internalLchaininfree}
Let $\pv V$ be a pseudovariety of monoids such that $\pv V=\pv
A\malce \pv V$ and let $X$ be a finite set. Then, for each
$\gamma\in \FP V X$, the submonoid $\stab {\gamma}$ is an internal
$\L$-chain.
\end{Thm}
\begin{proof}
Since $X$ is finite, we may write $\FP V X = \ilim_{n\in \mathbb N}
M_n$ where the $M_n$ are finite $X$-generated monoids in \pv V.
Let $\pi_n:\FP V
X\to M_n$ be the canonical projection.  Then, for $\gamma\in
\FP V X$, we have $\stab {\gamma} = \ilim_{n\in \mathbb N} \stab
{\gamma\pi_n}$, by Lemma~\ref{canonical}.  Let $\delta,\sigma\in \stab{\gamma}$ and consider
$M_n$. Since $\pv V=\pv A\malce \pv V$, we have that $\hatexp
2 {M_n}\in \pv V$.   Then $[\delta]_{\hatexp 2 {M_n}},
[\sigma]_{\hatexp 2 {M_n}}\in \stab{[\gamma]_{\hatexp 2
{M_n}}}$ and so Lemma~\ref{criticallemmaforpairs} implies
that $[\delta]_{M_n},[\sigma]_{M_n}$ are comparable in
the $\leq_{\L}$-order on $\stab {[\gamma]_{M_n}}$.  By going
to a subsequence we may assume without loss of generality  that, say,
$[\delta]_{M_n}\leq_{\L} [\sigma]_{M_n}$ in $\stab
{[\gamma]_{M_n}}$ for all $n$.  It then follows that $\delta
\leq_{\L}\sigma$ in $\stab {\gamma}$ (c.f.~\cite[Theorem
5.6.1]{Almeida:book} or~\cite[Proposition 9.1]{RhodesStein}).  Hence
$\stab {\gamma}$ is an
internal $\L$-chain.
\end{proof}

In~\cite[Corollary 14.5]{RhodesStein} it was shown that, for $\alpha\in
\wh{X^*}$, $\stab \alpha$ is an $\R$-trivial band. We can now refine
this result.

\begin{Cor}\label{Rtrivband}
Let $X$ be a finite set and $\alpha\in \wh{X^*}$.  Then $\stab
\alpha$ is an $\L$-chain of idempotents.  In particular, it is an
$\R$-trivial band.
\end{Cor}

\subsection{A discontinuous homomorphism}

The next lemma is the principal advantage obtained by our profinite
approach over Henckell's approach~\cite{Henckellstable}.

\begin{Lemma}\label{canonicalrelmorphiiseq}
Let $\pv V$ be a pseudovariety of monoids such that $\pv V=\pv
A\malce \pv V$.  Let $X$ be a finite set, $M$ an $X$-generated
finite monoid and $\rho_{\pv V}:M\to \FP V X$ the canonical
relational morphism.  Then \[\gamma\rho_{\pv V}\inv\sigma\rho_{\pv
V}\inv = (\gamma\sigma)\rho_{\pv V}\inv\] for all $\gamma,\sigma\in \FP
V X$.
\end{Lemma}
\begin{proof}
The inclusion $\gamma\rho_{\pv V}\inv\sigma\rho_{\pv V}\inv
\subseteq (\gamma\sigma)\rho_{\pv V}\inv$ is true for any relational
morphism. Since $1\in \varepsilon\rho_{\pv V}\inv$, the reverse
inclusion is trivial if either $\sigma$ or $\gamma$ is
$\varepsilon$, so assume $\sigma\neq\varepsilon\neq \gamma$.   Let
$m\in (\gamma\sigma)\rho_{\pv V}$. Then there exists a sequence of
words $w_n\in X^+$ such that $w_n\to \gamma\sigma$ and $[w_n]_M =
m$. Since $X$ is a finite set, we can write $\FP V X=\ilim_{n\in
\mathbb N} M_n$ where the $M_n$ are $X$-generated monoids from $\pv
V$. Again, $\pv V$ is closed under the expansion $N\mapsto \hatexp 2
{N}$.  By going to a subsequence, we may assume that $[w_n]_{\hatexp
2 {M_n}} = [\gamma\sigma]_{\hatexp 2 {M_n}}$, all $n$. Similarly, we
can find sequences $u_n,v_n\in X^+$ such that $u_n\to \gamma$,
$v_n\to \sigma$ and $[u_n]_{\hatexp 2 {M_n}} = [\gamma]_{\hatexp 2
{M_n}}$, $[v_n]_{\hatexp 2 {M_n}} = [\sigma]_{\hatexp 2 {M_n}}$.
Hence $[u_nv_n]_{\hatexp 2 {M_n}}=[\gamma\sigma]_{\hatexp 2
{M_n}}=[w_n]_{\hatexp 2 {M_n}}$ and so, by
Proposition~\ref{congruencehat}, $w_n = c_ns_n$ with $[c_n]_{M_n} =
[u_n]_{M_n} = [\gamma]_{M_n}$ and $[s_n]_{M_n} = [v_n]_{M_n} =
[\sigma]_{M_n}$.  Thus $c_n\to \gamma$ and $s_n\to \sigma$.  Since
$M$ is finite, by going to a subsequence, we may assume that
$[c_n]_M$ and $[s_n]_M$ are constant, say $[c_n]_M = m_1$ and
$[s_n]_M=m_2$. Then $m_1\in \gamma\rho_{\pv V}\inv$, $m_2\in
\sigma\rho_{\pv V}\inv$ and $m_1m_2 = [c_ns_n]_M = [w_n]_M = m$.  So
$m\in \gamma\rho_{\pv V}\inv\sigma\rho_{\pv V}\inv$.  This
establishes $(\gamma\sigma)\rho_{\pv V}\inv\subseteq \gamma\rho_{\pv
V}\inv\sigma\rho_{\pv V}\inv$ and completes the proof of the lemma.
\end{proof}

Let us reformulate the above result into our critical lemma.

\begin{Lemma}\label{Homolemma}
Let $M$ be an $X$-generated finite monoid, let \pv V be a pseudovariety
such that $\pv V=\pv A\malce \pv V$ and let $\rho_{\pv V}:M\to \FP V
X$ be the canonical relational morphism.  Then the map $f_{\pv V}:\FP V
X\to \pl V M$ defined by $\gamma f_{\pv V}= \gamma\rho_{\pv V}\inv$ is
a monoid homomorphism.
\end{Lemma}
\begin{proof}
Theorem~\ref{profinitewitness} shows that $f_{\pv V}$ is well
defined.  Lemma~\ref{canonicalrelmorphiiseq} shows that $f_{\pv V}$
is a semigroup homomorphism.  Since $\varepsilon$ is an isolated
point of $\FP V X$ (as the congruence class of $\varepsilon$ is
trivial in $\hatexp 2 M$ for any $M\in \pv V$), we conclude
$\varepsilon f_{\pv V}=\varepsilon\rho_{\pv V}\inv =\{1\}$ and thus
$f_{\pv V}$ is a monoid homomorphism.
\end{proof}

We remark that $f_{\pv V}$ is not necessarily continuous.  For
instance, if $\gamma f_{\pv A} = Z$, then $\gamma^{\omega} f_{\pv
A}\supseteq \bigcup_{n\in \mathbb N}Z^{\omega}Z^n$, which can be
strictly bigger than $Z^{\omega}$.  This discontinuity is what
underlies the analysis of aperiodic pointlike
sets~\cite{Henckell,ourAPL}.

\subsection{Idempotent pointlikes}
As a warm-up we prove the result of Henckell~\cite{Henckellidem}
relating \pv V-idempotent pointlikes with idempotent \pv
V-pointlikes.

\begin{Thm}[Henckell]
Let $\pv V$ be a pseudovariety of monoids such that $\pv A\malce \pv V=\pv V$
and let $M$ be a finite monoid.  Then
the maximal $\pv V$-idempotent pointlikes of $M$ are precisely the
maximal idempotents of $\pl V M$.
\end{Thm}
\begin{proof}
We already observed that idempotents of $\pl V M$ are \pv V-idempotent
pointlike.  Conversely, suppose that $Z$ is a maximal \pv V-idempotent
pointlike subset of $M$.  Let $X$ be a
finite generating set for $M$ and let
\mbox{$\rho_{\pv V}:M\to \FP V X$} and $f_{\pv V}:\FP V X\to \pl V M$ be as
per Lemma~\ref{Homolemma}. By
maximality and Theorem~\ref{profinitewitness}, we must have that $Z =
e\rho_{\pv V}\inv = ef_{\pv V}$ for some idempotent $e\in \FP V
X$. Since $f_{\pv V}$ is a homomorphism, $Z\in \pl V A$ is idempotent.
\end{proof}

Since Henckell proved~\cite{Henckell,ourAPL} that $\pv A$-pointlikes are
decidable, we have the following corollaries.

\begin{Cor}
If $\pv V = \pv A\malce \pv V$ and $\pv V$-pointlikes are decidable,
then $\pv V$-idempotent pointlikes are decidable.  In particular,
$\pv A$-idempotent pointlikes are decidable.
\end{Cor}

\begin{Cor}
If \pv V has decidable membership, then do does $\pv V\malce \pv A$.
\end{Cor}

\section{\pv M-stable pairs}
Let \pv M be the pseudovariety of all finite monoids.  As a second
warm-up exercise we characterize the $\pv M$-stable pairs.  This
partially answers a question raised in~\cite{Almeidahyp}. By
considering the identity map, we see that an \pv M-stable pair of a
monoid $M$ must be of the form $(\{x\},N)$ where $N\leq \stab x$.
Corollary~\ref{Rtrivband} suggests that $\L$-chains of idempotents
should play a role.  The next lemma describes what kind of monoid
you can obtain by such a chain.

\begin{Lemma}\label{anotherlrbresult}
Suppose that $e_1\geq _{\L} e_2\geq_{\L}\cdots \geq_{\L} e_n$ is an
$\L$-chain of idempotents in a monoid $M$.  Then $\langle
e_1,\ldots,e_n\rangle$ is an $\R$-trivial band.
\end{Lemma}
\begin{proof}
Set $N=\langle e_1,\ldots,e_n\rangle$.  First we observe that
$e_ie_j = e_i$ if $i\geq j$. Thus, each element of $T$ can be
written in the form $e_{i_1}e_{i_2}\cdots e_{i_m}$ where the indices
are increasing: $i_1<i_2<\cdots<i_m$. Clearly then one has
\[(e_{i_1}e_{i_2}\cdots e_{i_m})^2= e_{i_1}e_{i_2}\ldots e_{i_m},\]
since $i_m\geq i_j$ for all $1\leq j\leq m$. Let
$s=e_{i_1}e_{i_2}\cdots e_{i_m}$ and $t=e_{j_1}e_{j_2}\cdots
e_{j_{\ell}}$.  Then $st = s$ if $i_m\geq j_{\ell}$, or else $st=
se_{j_k}\cdots e_{j_\ell}$, where $k$ is the smallest index such
that $i_m< j_k$. In the first case, clearly $sts=s^2=s=st$, while in
the latter case we have $j_{\ell}\geq i_r$ for all $r$ and so
$e_{j_{\ell}}s = e_{j_{\ell}}$, from which we conclude that
$sts=st$. This proves that $N$ is an $\R$-trivial band.
\end{proof}

\begin{Thm}
Let $M$ be a finite monoid.  Then $(\{y\},N)$ is an \pv M-stable
pair of $M$ if and only if there there is an $\L$-chain $Y$ of
idempotents in $\stab y$ such that $N\leq \langle Y\rangle$.
\end{Thm}
\begin{proof}
Suppose that $Y\subseteq \stab y$ is an $\L$-chain of idempotents.
Without loss of generality, we may assume that $N=\langle Y\rangle$
(since stable pairs are downwards closed). Since every monoid
belongs to $\pv M$, it clearly suffices to show that if
$\p:S\twoheadrightarrow M$ is an onto homomorphism, then there
exists $y'\in S$ such that $y'\p = y$ and $N\leq \stab {y'}\p$.

Choose $\til y\in S$ with $\til y\p = y$.  Next, suppose $Y =e_1\geq
_{\L} e_2\geq_{\L}\cdots \geq_{\L} e_n$ and choose an idempotent
$f_1\in S$ with $f_1\p = e_1$.  Assume inductively that, for
\mbox{$1\leq i<n$} we have found $f_1\geq_{\L}\cdots\geq_{\L} f_i$
in $S$ with $f_j\p = e_j$, for \mbox{$1\leq j\leq i$}.  Then
$e_{i+1}\in Me_i\subseteq (Sf_i)\p$. So there exists an idempotent
$f_{i+1}$ of $Sf_i$ with $f_{i+1}\p = e_{i+1}$. This completes the
induction. Set $Y' = \{f_1,\ldots,f_n\}$. Clearly $N' = \langle
Y'\rangle$ maps onto $N$ via $\p$.  Also $N'$ is an $\R$-trivial
band by Lemma~\ref{anotherlrbresult}.  In particular, if $s$ belongs
to the minimal ideal of $N'$ and $t\in N'$, then $st\R s$ and hence,
since $N'$ is $\R$-trivial, $st=s$. Thus $N'\leq \stab s\leq
\stab{\til ys}$.  But $(\til ys)\p = \til y\p s\p =ys\p\in
yN=\{y\}$.  This shows that $(\{y\},N)$ is an $\pv M$-stable pair.

For the converse, choose a generating set $X$ for $M$.  Let
$\pi:\wh{X^*}\to M$ be the canonical projection.  Then $\pi\inv$ is
the canonical relational morphism $\rho_{\pv M}:M\to \wh{X^*}$.  So
Theorem~\ref{profinitewitness} shows there exists $\alpha\in
\wh{X^*}$ with $\alpha\pi = y$ and $N\leq \stab {\alpha}\pi$.
Corollary~\ref{Rtrivband} yields $\stab {\alpha}$ is an $\L$-chain
of idempotents, from which the result easily follows.
\end{proof}

\section{\pv A-stable pairs}
The situation for $\pv A$-stable pairs is more complicated since we
no longer have that the stabilizers in $\FP A X$ must be bands.  Let
us recall some terminology from~\cite{lowerbounds1,lowerbounds2}
(see also~\cite{qtheor}).  Let $\pv {ER}$ be the pseudovariety of
monoids whose idempotent-generated submonoids are $\R$-trivial. It
is well-known that $M\in \pv {ER}$ if and only if each regular
$\R$-class of $M$ contains a unique idempotent~\cite{qtheor}.

\begin{Prop}\label{ERaperiodic}
Let $M\in \pv {ER}\cap \pv A$.  Then, for any $x$ in the minimal
ideal of $M$, one has $\stab x= M$.
\end{Prop}
\begin{proof}
Since $M\in \pv {ER}$, the minimal ideal $I$ of $M$ contains a
unique $\L$-class.  If $x\in I$ and $m\in M$, then $xm\R x$ by
stability of finite semigroups and $xm\L x$ since $I$ has a unique
$\L$-class.  Since $M$ is aperiodic, $xm=x$.
\end{proof}

A monoid $M$ is said to be \emph{absolute Type
I}~\cite{lowerbounds1,lowerbounds2,qtheor,HMPR} if it can be
generated by a chain of its $\L$-classes.  In particular, an
internal $\L$-chain is absolute Type I.  The facts contained in our
next proposition are from~\cite{lowerbounds1}; see~\cite{qtheor} for
proofs.

\begin{Prop}\label{T1}
{}\ \begin{enumerate}
\item An aperiodic absolute Type I-monoid belongs to $\pv {ER}$.
\item If $\p:M\twoheadrightarrow N$ is an onto homomorphism and $M$
is absolute Type I, then $N$ is absolute Type I.
\item If $\p:M\twoheadrightarrow N$ is an onto homomorphism and $N$
is absolute Type I, then there is an absolute Type I-submonoid $M'\leq
M$ with $M'\p = N$.
\end{enumerate}
\end{Prop}

We now present a sufficient condition for  $(Y,N)$ to be an
$\pv A$-stable pair for a monoid $M$.
\begin{Prop}\label{suff}
Let $M$ be a finite monoid.  Suppose that $Y\in \pl A M$ and $W\leq
\pl A M$ is a submonoid which is an internal $\L$-chain such that:
\begin{enumerate}
\item $\bigcup W = N$;
\item $W\leq \stab {Y}$.
\end{enumerate}
Then $(Y,N)$ is an $\pv A$-stable pair.
\end{Prop}
\begin{proof}
Let $\p:M\to A$ with $A\in \pv A$ be a relational morphism.  Factor
$\p=\alpha\inv \beta$ where $\alpha:R\twoheadrightarrow M$ is an
onto homomorphism, $\beta:R\to A$ is a homomorphism and $R$ is
finite.  Let $\alpha_*:\pl A R\twoheadrightarrow \pl A M$ and
$\beta_*:\pl A R\to \pl A A = A$ be the induced maps from
Proposition~\ref{liftandpushpoint}.  We shall use several times that
if  $X\beta_*=x$, then $X\subseteq x\beta\inv$.  Since $W$ is
absolute Type I, we can find, by Proposition~\ref{T1}, an absolute
Type I submonoid $W'\leq \pl A R$ with $W'\alpha_* = W$. Then
$W''=W'\beta_*$ is absolute Type I and hence belongs to $\pv
{ER}\cap \pv A$ (again by Proposition~\ref{T1}).  Choose $a\in A$
with $Y\subseteq a\pinv$ and choose $z''$ from the minimal ideal of
$W''$.  By definition of $W''$, there exists $Z'\in W'$ with
$Z'\beta_*= z''$. Setting $Z= Z'\alpha_*\in W$, we have $Z =
Z'\alpha \subseteq z''\beta\inv \alpha = z''\pinv$ and so, as $W\leq
\stab {Y}$,
\begin{equation}\label{firstsuff}
Y=YZ\subseteq a\pinv z''\pinv \subseteq (az'')\pinv.
\end{equation}

  Now
Proposition~\ref{ERaperiodic} shows that $W''\subseteq \stab
{z''}\subseteq \stab {az''}$. So we are left with showing that
$N\leq W''\pinv$.  Let $U\in W$.  Then we can find $U'\in W'$ such
that $U'\alpha_*=U$. Set $U'\beta_*=u''\in W''$.  Then $U = U'\alpha
\subseteq u''\beta\inv\alpha = u''\pinv$. Thus $U\subseteq W''\pinv$
and so we may conclude
\begin{equation}\label{secsuff}
N = \bigcup W \subseteq W''\pinv \subseteq \stab {az''}\pinv.
\end{equation}
Combining \eqref{firstsuff} and \eqref{secsuff} yields that $(Y,N)$
is an \pv A-stable pair.
\end{proof}

We now prove the converse for maximal \pv A-stable pairs;
Henckell proves an apparently stronger formulation
in~\cite{Henckellstable}.

\begin{Thm}\label{main}
Suppose that $M$ is a finite monoid.  Then the maximal \pv A-stable
pairs of $M$ are the maximal pairs $(Y,N)$ such that $Y\in \pl A M$
and there exists a submonoid \mbox{$W\leq \pl A M$} with $W$  an
internal $\L$-chain and:
\begin{enumerate}
\item $\bigcup W = N$;
\item $W\leq \stab Y$.
\end{enumerate}
\end{Thm}
\begin{proof}
By Proposition~\ref{suff} any such pair $(Y,N)$ is \pv A-stable.
Conversely, suppose that $(Y,N)$ is a maximal \pv A-stable pair for
$M$.  Choose a finite generating set $X$ for $M$ and let $\rho_{\pv
A}:M\to \FP A X$ be the canonical relational morphism. Let $f_{\pv
A}:\FP A X\to \pl A M$ be the homomorphism from
Lemma~\ref{Homolemma}; so $f_{\pv A} = \rho_{\pv A}\inv$. Maximality
and Theorem~\ref{profinitewitness} implies there exists $\gamma\in
\FP A X$ such that $Y= \gamma\rho_{\pv A}\inv=\gamma f_{\pv
  A}$ and \mbox{$N = \stab {\gamma}\rho_{\pv A}\inv$}.   By
Theorem~\ref{internalLchaininfree}, $\stab {\gamma}$ is an internal
$\L$-chain.   Then we see that \[W=\stab {\gamma}f_{\pv A}\leq \stab
{\gamma f_{\pv A}} = \stab Y\] is a submonoid of $\pl A M$ and an
internal $\L$-chain. Moreover, we have \[\bigcup W = \bigcup
_{\beta\in\stab{\gamma}} \beta f_{\pv A} =\bigcup
_{\beta\in\stab{\gamma}} \beta \rho_{\pv A}\inv = \stab
{\gamma}\rho_{\pv A}\inv =N.\] This completes the proof of the
theorem.
\end{proof}

Since $\pl A M$ is computable~\cite{Henckell,ourAPL},
Theorem~\ref{main} admits as corollaries:

\begin{Cor}
Stable pairs are decidable for $\pv A$.  Equivalently, \pv
A-in\-ev\-i\-ta\-bil\-i\-ty is decidable for labellings of graphs
with a single vertex, with singletons on the edges.
\end{Cor}

\begin{Cor}
If \pv V is a local pseudovariety with decidable membership, then
$\pv V\ast \pv A$ is decidable.
\end{Cor}

\section{\pv A-triples}
To compute the Krohn-Rhodes complexity of a monoid, we shall need
some other notions, related to those we have been considering.

\begin{Def}[\pv V-triple]
Let us call a triple $(A,B,C)$ of subsets of a finite monoid $M$ a
$\pv V$-triple if, for all relational morphisms $\p:M\to N$ with
$N\in \pv V$, there exist $a,b,c\in N$ such that $A\subseteq
a\pinv$, $B\subseteq c\pinv$, $C\subseteq c\pinv$ and $abc=ab$.
\end{Def}

Equivalently, $(A,B,C)$ is a \pv V-triple if and only if the graph
with two vertices $v_1, v_2$, an edge $e_1$ from $v_1$ to $v_2$ and
a loop $e_2$ from $v_2$ to $v_2$ with $v_1$ labelled by $A$, $e_1$
by $B$, $v_2$ by $AB$ and $e_2$ by $C$ is $\pv V$-inevitable. Thus
an analogue of Theorem~\ref{profinitewitness} holds for $\pv V$-triples.

We are particularly interested in $\pv A$-triples and so we begin by
investigating solutions to equations of the from $xyz=xy$ in $\FP A
X$. It turns out that the Henckell-Sch\"utzenberger expansion allows
one to treat equations over $\FP A X$ in a similar way to equations
over free monoids.

\begin{Prop}\label{tripleinfreeproA}
Let $X$ be a finite set and let $\alpha,\beta,\gamma \in \FP A X$.
Then $\alpha\beta\gamma = \alpha\beta$ if and only if one of the
following three situations occur:
\begin{enumerate}
\item $\beta\gamma = \beta$;
\item there exists $\tau\in \FP A X$ such that $\alpha\beta\tau = \alpha$
and $\gamma = \tau\beta$;
\item there exist $\sigma,\tau\in \FP A X$ and $i\geq 1$ such that $\alpha =
\alpha \tau\sigma$, $\beta = (\tau\sigma)^i\tau$ and $\gamma =
\sigma\tau$.
\end{enumerate}
\end{Prop}
\begin{proof}
Clearly any of (1), (2) or (3) implies
$\alpha\beta\gamma=\alpha\beta$.  For the converse, if $\alpha$ or
$\gamma$ are $\varepsilon$, we are in case (1). If $\beta
=\varepsilon$, then we are in case (2) with $\tau=\gamma$. So we may
assume that none of $\alpha$, $\beta$ and $\gamma$ are
$\varepsilon$.   Since $X$ is finite, we may write $\FP A X =
\ilim_{n\in \mathbb N} M_n$ with the $M_n$ finite $X$-generated
aperiodic monoids. Moreover, $\hatexp 2 {M_n}\in \pv A$ for all $n$.
Choose sequences of words $a_n,b_n,c_n$ from $X^+$ such that $a_n\to
\alpha$, $b_n\to \beta$ and $c_n\to \gamma$. By passing to
subsequences, we may assume that $[a_n]_{\hatexp 2 {M_n}} =
[\alpha]_{\hatexp 2 {M_n}}$, $[b_n]_{\hatexp 2 {M_n}} =
[\beta]_{\hatexp 2 {M_n}}$ and $[c_n]_{\hatexp 2 {M_n}} =
[\gamma]_{\hatexp 2 {M_n}}$, for all $n$.

Then, for each $n$, we have the equality $[a_nb_nc_n]_{\hatexp 2
{M_n}} = [a_nb_n]_{\hatexp 2 {M_n}}$.
Proposition~\ref{congruencehat} says that $a_nb_nc_n = a_n'b_n'$
with $[a_n']_{M_n} = [a_n]_{M_n} = [\alpha]_{M_n}$ and $[b_n']_{M_n}
= [b_n]_{M_n} = [\beta]_{M_n}$. In particular, $a_n'\to \alpha$ and
$b_n'\to \beta$.

For each $n$, there are three cases: $|a_n'|\leq |a_n|$, $|a_n'|\geq
|a_nb_n|$ and finally \mbox{$|a_n|<|a_n'|<|a_nb_n|$}. By passing to
a subsequence, we may assume that the same case occurs for all $n$.

Suppose that $|a_n'|\leq |a_n|$ for all $n$.  Then, for each $n$,
there exists $t_n\in X^*$ so that $a_n = a_n't_n$, $b_n' =
t_nb_nc_n$.  By passing to a subsequence, we may assume that $t_n\to
\tau \in \FP A X$.  Then $\beta = \tau \beta \gamma$ and so we have
$\beta = \tau^{\omega}\beta \gamma^{\omega}$.  Thus $\beta \gamma =
\tau^{\omega}\beta \gamma^{\omega}\gamma = \tau^{\omega}\beta
\gamma^{\omega}=\beta$ and we are in case (1).

Next suppose that $|a_n'|\geq |a_nb_n|$ for all $n$.  Then, for each
$n$, we can find $t_n\in X^*$ such that $a_n' = a_nb_nt_n$ and $c_n
= t_nb_n'$.  By passing to a subsequence, we may assume $t_n\to
\tau\in \FP A X$.  Then $\alpha = \alpha\beta \tau$ and $\gamma =
\tau \beta$, and so we are in case (2).

Finally, suppose $|a_n|<|a_n'|<|a_nb_n|$ for all $n$.  Then, for
each $n$, we can find $p_n,t_n\in X^*$ such that $a_n' = a_np_n$,
$b_n' = t_nc_n$ and $b_n = p_nt_n$.  By extracting a subsequence, we
may assume that $p_n\to \pi$ and $t_n\to \tau_1$ in $\FP A X$.  Then
we have in $\FP A X$ the equalities $\tau_1\gamma = \beta =
\pi\tau_1$ and
\begin{equation}\label{tripleseq1}
\alpha = \alpha\pi.
\end{equation}

Define $\tau_0=\beta$ and suppose inductively that we have found
\mbox{$\tau_i\in \FP A X$}, for $i\geq 1$, such that $\tau_i\gamma =
\tau_{i-1}= \pi\tau_i$. Notice that a simple induction yields
\begin{equation}\label{tripleseq2}
\tau_i\gamma^i=\beta=\pi^i\tau_i
\end{equation}
  Then we can choose sequences of
words $P_n,T_n$ such that $P_n\to \pi$, \mbox{$T_n\to \tau_i$} and
$[P_n]_{\hatexp 2 {M_n}} = [\pi]_{\hatexp 2 {M_n}}$, $[T_n]_{\hatexp
2 {M_n}} = [\tau_i]_{\hatexp 2 {M_n}}$, all $n$. Then, for each $n$,
we have
\[[T_nc_n]_{\hatexp 2 {M_n}}=[\tau_i\gamma]_{\hatexp 2 {M_n}}=
[\pi\tau_i]_{\hatexp 2 {M_n}}= [P_nT_n]_{\hatexp 2 {M_n}}.\]
Proposition~\ref{congruencehat} then shows  $T_nc_n = P_n'T_n'$
where \[[P_n']_{M_n} = [P_n]_{M_n}=[\pi]_{M_n},\quad [T_n']_{M_n} =
[T_n]_{M_n}=[\tau_i]_{M_n}.\]  In particular, we have $P_n'\to \pi$
and $T_n'\to \tau_i$.

For any $n$, there are two cases: $|P_n'|< |T_n|$ and $|P_n'|\geq
|T_n|$.  By passing to a subsequence, we may assume that the same
case applies for all $n$.  Suppose first that $|P_n'|<|T_n|$ for all
$n$. Then we can find $R_n\in X^*$ so that $T_n = P_n'R_n$ and $T_n'
= R_nc_n$. Extracting a subsequence, we may assume that $R_n$
converges to some $\tau_{i+1}$ in $\FP A X$. Then
$\pi\tau_{i+1}=\tau_i= \tau_{i+1}\gamma$, allowing us to continue
the induction.

Next assume that $|P_n'|\geq |T_n|$ for all $n$.  Then $P_n' =
T_nS_n$ and $c_n = S_nT_n'$ for some $S_n\in X^*$.  By passing to a
subsequence, we may assume that $S_n\to \sigma$ in $\FP A X$.  Then
$\pi = \tau_i \sigma$ and $\gamma = \sigma\tau_i$.  Therefore, by
\eqref{tripleseq1} and \eqref{tripleseq2}, we have the equalities
\[\alpha = \alpha\tau_i\sigma,\quad \beta = \pi^i\tau_i =
(\tau_i\sigma)^i\tau_i,\quad \gamma = \sigma\tau_i,\] and so we are
in case (3) and may stop.

Hence, either one of cases (1), (2) or (3) arises, or we can find an
infinite sequence $\{\tau_i\}$ of elements of $\FP A X$ with $\beta
= \tau_i\gamma^i$.  By passing to a subsequence, we may assume that
$\tau_i\to \tau\in \FP A X$.  Since $\lim_{i\to\infty}\gamma^i=
\gamma^{\omega}$, we obtain $\beta = \tau\gamma^{\omega}$ and hence
$\beta\gamma = \tau\gamma^{\omega}\gamma = \tau\gamma^{\omega} =
\beta$, so we are again in case (1). This completes the proof.
\end{proof}

\begin{Cor}
Let $M$ be a finite monoid.  Then the maximal $\pv A$-triples are
the maximal triples $(A,B,C)\in \pl A M^3$ such that one of the
following occurs:
\begin{enumerate}
\item $BC=B$;
\item there exists $T\in \pl A M$ such that $ABT=A$ and $C = TB$;
\item there exist $S,T\in \pl A M$ and $i\geq 1$ such that $A=ATS$,
$B = (TS)^iT$ and $C = ST$.
\end{enumerate}
In particular, $\pv A$-triples are decidable.
\end{Cor}
\begin{proof}
First we show that if $(A,B,C)\in \pl A M^3$ satisfies any of
(1)--(3), then it is an \pv A-triple.  Let $\p:M\to N$ with $N\in
\pv A$ be a relational morphism.

Suppose that (1) holds.  Choose $a,b,c\in N$ with $A\subseteq
a\pinv$, $B\subseteq b\pinv$ and $C\subseteq c\pinv$.  Then
$B\subseteq bc^{\omega}\pinv$ and $a(bc^{\omega})c =
a(bc^{\omega})$. Thus $(A,B,C)$ is an $\pv A$-triple.

Next assume that (2) holds. Choose $a,b,t\in N$ with $A\subseteq
a\pinv$, $B\subseteq b\pinv$ and $T\subseteq t\pinv$. Then we have
$A\subseteq a(bt)^{\omega}\pinv$, $C\subseteq tb\pinv$ and the
equality $[a(bt)^{\omega}]b(tb) = [a(bt)^{\omega}]b$, and so
$(A,B,C)$ is an \pv A-triple.

Finally, assume that (3) holds. Choose $a,s,t\in N$ with $A\subseteq
a\pinv$, \mbox{$S\subseteq s\pinv$} and $T\subseteq t\pinv$.  Then
we have $A\subseteq a(ts)^{\omega}\pinv$, $B\subseteq (ts)^it\pinv$
and \mbox{$C\subseteq st\pinv$}.  Moreover,
$[a(ts)^{\omega}][(ts)^it](st) = [a(ts)^{\omega}][(ts)^it]$.  So we
see that in all cases $(A,B,C)$ is an \pv A-triple.

Next suppose that $(A,B,C)$ is a maximal \pv A-triple.  Choose a
generating set $X$ for $M$ and let $\rho_{\pv A}:M\to \FP A X$ be
the canonical relational morphism.  Let $f_{\pv A}:\FP A X\to \pl A
M$ be the homomorphism from Lemma~\ref{Homolemma}.  By
Theorem~\ref{profinitewitness} and maximality, we can find
$\alpha,\beta,\gamma\in \FP A X$ such that
\begin{equation*}
A = \alpha\rho_{\pv A}\inv = \alpha f_{\pv A},\quad B =
\beta\rho_{\pv A}\inv = \beta f_{\pv A},\quad C= \gamma\rho_{\pv
A}\inv = \gamma f_{\pv A}.
\end{equation*}
and $\alpha\beta\gamma = \alpha\beta$.  We analyze the situation
according to the three cases of Proposition~\ref{tripleinfreeproA}.
If $\beta\gamma=\beta$, then \[BC=\beta f_{\pv A}\gamma f_{\pv A} =
(\beta\gamma)f_{\pv A} = \beta f_{\pv A}=B\] and we are in case (1).
If there exists $\tau\in \FP A X$ such that $\alpha\beta\tau
=\alpha$ and $\gamma = \tau\beta$, then setting $T=\tau f_{\pv A}$,
we have
\begin{align*}
A&= \alpha f_{\pv A} = (\alpha\beta\tau)f_{\pv A} = \alpha f_{\pv
A}\beta f_{\pv A}\tau f_{\pv A} = ABT\\
C &= \gamma f_{\pv A} = (\tau\beta)f_{\pv A} = \tau f_{\pv A}\beta
f_{\pv A} = TB,
\end{align*}
and so we are in case (2).

Finally, if there exist $\sigma,\tau\in \FP A X$ and $i\geq 1$ such
that $\alpha = \alpha \tau\sigma$, $\beta = (\tau\sigma)^i\tau$ and
$\gamma = \sigma\tau$, then setting $S=\sigma f_{\pv A}$ and $T=\tau
f_{\pv A}$, we have
\begin{align*}
A &= \alpha f_{\pv A} = (\alpha\tau\sigma)f_{\pv A}=\alpha f_{\pv
A}\tau f_{\pv A}\sigma f_{\pv A} = ATS\\
B &= \beta f_{\pv A} = ((\tau\sigma)^i\tau)f_{\pv A} =(\tau f_{\pv
A}\sigma f_{\pv A})^i\tau f_{\pv A} = (TS)^iT\\
C &= \gamma f_{\pv A} = (\sigma\tau)f_{\pv A} = \sigma f_{\pv A}\tau
f_{\pv A} = ST,
\end{align*}
and hence we are in case (3).  This completes the proof.
\end{proof}

\bibliographystyle{abbrv}
\bibliography{standard}

\end{document}